%Formatting macros
\baselineskip 24pt

\def\Q{{\bf {Q}}}
\def\K{{\bf {K}}}
\def\G{{\bf {G}}}

\def\Z{{\bf Z}}     
\def\R{{\bf R}}   \def\cK{{\bf B}}

\def\cK{{\cal {K}}}  
\def\cH{{\cal {H}}} \def\cL{{\cal {L}}} \def\cV{{\cal {V}}}

\def\cM{{\cal {M}}}

\def\a{{\bf {a}}}
\def\j{{\bf {j}}}
\def\r{{\bf {r}}}

\magnification=1200

\def\house#1{\setbox1=\hbox{$\,#1\,$}%
\dimen1=\ht1 \advance\dimen1 by 2pt \dimen2=\dp1 \advance\dimen2 by 2pt
\setbox1=\hbox{\vrule height\dimen1 depth\dimen2\box1\vrule}%
\setbox1=\vbox{\hrule\box1}%
\advance\dimen1 by .4pt \ht1=\dimen1
\advance\dimen2 by .4pt \dp1=\dimen2 \box1\relax}

\def\HH{{\rm H}}
    
  \def\de{\delta _{\K}}

   \def\al{{\alpha}}

  \def\de{{\delta}}

\def\cH{{\cal H}}

\def\sm{\smallskip} \def\ens{\enspace} \def\noi{\noindent}

\def\build#1_#2^#3{\mathrel{\mathop{\kern 0pt#1}\limits_{#2}^{#3}}}

\def\date {le\ {\the\day}\ \ifcase\month\or janvier
\or fevrier\or mars\or avril\or mai\or juin\or juillet\or
ao\^ut\or septembre\or octobre\or novembre
\or d\'ecembre\fi\ {\oldstyle\the\year}}

\font\fivegoth=eufm5 \font\sevengoth=eufm7 \font\tengoth=eufm10

\newfam\gothfam \scriptscriptfont\gothfam=\fivegoth
\textfont\gothfam=\tengoth \scriptfont\gothfam=\sevengoth
\def\goth{\fam\gothfam\tengoth}

\def\BA{{\goth B}}

\def\smallsquare{\vbox{\hrule\hbox{\vrule height 1 ex\kern 1
ex\vrule}\hrule}}
\def\cqfd{\hfill \smallsquare\vskip 3mm}

%%%%%%%%%%%%%%%%%%%%%%%%%%%%%%%%%%%%%%%%%%%%

\centerline{}

\vskip 4mm

\centerline{
\bf Zero--infinity laws in Diophantine approximation}

\vskip 8mm
\centerline{Y. B{\sevenrm UGEAUD}, \footnote{}{\rm 
} \
M. M. D{\sevenrm ODSON} \ \&  \ S. K{\sevenrm
RISTENSEN}\footnote{}{Research partially funded by a CNRS/Royal
Society exchange grant. SK is a William Gordon Seggie Brown
Fellow.}}

\vskip 6mm

{\narrower\narrower
\vskip 12mm

\proclaim Abstract. {It is shown that for any translation invariant
outer measure $\cM$, the $\cM$-measure of the intersection of any
subset of \ $\R^n$ that is invariant under rational translations and
which does not have full Lebesgue measure with an the closure of an
open set of positive measure cannot be positive and finite.  Analogues
for $p$-adic fields and fields of formal power series over a finite
field are established. The results are applied to some problems in
metric Diophantine approximation.}

}

\vskip 6mm
\vskip 8mm

\centerline{\bf 1. Introduction}

\vskip 6mm

`Zero--infinity' laws for the measure of a set are a natural analogue
of the more familiar `zero--one' laws of probablility theory.  They
arise in the setting of the real line or with measures other than
Lebesgue. Unlike a `zero--one' law, a `zero--infinity' law does not
imply that a set is null or full (or has null complement), so in this
respect, `zero--infinity' laws are weaker than the familiar
`zero--one' laws. This of course is a consequence of the special
properties of the Lebesgue measure.

For clarity, we will say that a set satisfies a {\it full}
`zero--infinity' law with respect to some measure if either the set or
its complement is null. It follows straightforwardly from the Lebesgue
density theorem (see {\it e.g.}\ [21], Lemma 7, p. 21) that such a
full `zero--infinity' law holds for Lebesgue measurable real sets
invariant under rational translations.

\proclaim Theorem A. {Let $E\subseteq\R$ be a Lebesgue measurable
set such that $\xi + p/q$ is in $E$ for any $\xi$ in
$E$ and any rational number $p/q$. Then either $E$ or its
complementary set has Lebesgue measure zero.}

One may ask whether Theorem A can be extended to other measures than
Lebesgue, and in particular to Hausdorff measures, for which there is
no analogue of the Lebesgue density theorem. Clearly, one cannot hope
for a complete analogue. Indeed, one may easily construct examples of
subsets of an interval where both the set and its complement has
infinite Hausdorff measure.  In view of this, we can ask that the set
$E$ be -- if not full -- at least `thick', {\it i.e.}, roughly
speaking for any `reasonable' set $F$, $\cM(E\cap F)\ge c\cM(F)$ for
some $c>0$.  This means that $E$ is sufficiently `spread out' to
ensure that the measure of the intersection with $F$ is always
comparable to the measure of $F$ (which in most cases will be
infinite). In particular, one can hope show that such a set must have
Hausdorff measure either zero or infinity for all dimension functions.

The first contribution to this problem is due to Jarn\'\i k [10, 11]
(see also [7], since these papers are not readily available).  He
showed that for any dimension function $f$ and any real number $\tau >
2$, the Hausdorff $f$-measure of the intersection of the unit interval
$[0,1]$ with
$$
\cV(\tau)  = \biggl\{\xi \in \R :  
\biggl| \xi - {p \over q} \biggr| < { K(\xi) \over q^\tau} \ens \hbox{ 
for infinitely many  }  {p \over q} \in \Q,
\hbox{ for some } K(\xi)>0 \biggr\},
$$ 
the set  of real numbers approximable by rational numbers to order
$\tau$,  is either $0$ or $+ \infty$. This set is invariant under
rational translations and is closely related to the set 
$$
\cK(\tau)  := \biggl\{\xi \in \R : 
\biggl| \xi - {p \over q} \biggr| < { 1 \over q^\tau} \ens \hbox{
for infinitely many  }  {p \over q} \in \Q  \biggr\},
$$ 
and indeed for any $\varepsilon>0$, we have $ \cV(\tau+\varepsilon)\subset
\cK(\tau)\subset \cV(\tau).  $ A little earlier Jarn\'\i k had shown
in~[9] that $\cH^s(\cK(\tau))$, the Hausdorff $f$-measure for $f(x) =
x^s$, was infinite when $s<2/\tau$ and vanished when $\tau>2/\tau$,
whence the Hausdorff dimension of $\cK(\tau)$ is given by
$$
\dim \cK(\tau)=2/\tau \ \hbox{ when } \tau>2
$$ 
($\cK(\tau)=\R$ when $\tau \le 2$).  It is not clear whether
$\cK(\tau)$ is invariant under rational translates but in his later
more general paper of 1931~[12] on the Hausdorff $s$-measure analogue
of Khintchine's theorem for simultaneous Diophantine approximation,
Jarn\'\i k proved that $\cH^{2/\tau}(\cK(\tau) \cap [0,1])= + \infty$.
(See {\it e.g.}\ Falconer [8] or Rogers [18] for background on the
theory of Hausdorff measure and dimension.) This shows that $\cV(\tau)
\cap [0,1]$ and $\cK(\tau) \cap [0,1]$ are not $s$-sets (see
Chapters~2--4 of~[8] for further details of $s$-sets), unlike for
instance the usual ternary Cantor set (whose Hausdorff $s$-measure at
the critical exponent $\log 2/\log 3$ is equal to $1$). The same
conclusion also holds for sets of real numbers well approximable by
algebraic numbers of bounded degree (see~[4]), the set $\BA$ of badly
approximable numbers (this set is invariant under rational translates,
see~[7]) and many other sets arising in Diophantine approximation.

It is thus natural to ask whether there exists a dimension function
$f$ for which the sets mentioned above and others intersected with a
non-trivial closed interval have $\cH^f$-measure that is always
positive and finite. For certain limsup sets such as $\cK(\tau)$, very
general zero-infinity statements for Hausdorff measures are known due
to Beresnevich, Dickinson and Velani [3], providing a negative
answer in the case of these sets.  Recently, Olsen [17] gave a
negative answer to this question for a different family of sets, under
the assumption that $f$ is strictly concave. The purpose of the
present paper is to remove this assumption by showing that Jarn\'\i
k's arguments actually imply that for any set $E$ in $\R^n$ invariant
under certain rational translations (such as the sets $\cV(\tau)$,
$\BA$ and the set $\cL$ of Liouville numbers), for any translation
invariant outer measure $\cM$, and for any set $F \subseteq \R$ which
is the closure of an open set with positive Lebesgue measure, either
$\cM(E \cap F) = 0$ or $\cM(E \cap F)= + \infty$. In particular, for
any dimension function $f$, we have $\cH^f (E \cap F) = 0$ or $+
\infty$.  Also, our general method provides a simple proof of Theorem
A.

Similar questions may be asked when the reals are replaced by a
$p$-adic field or a field of formal Laurent series with coefficients
from a finite field $\K$. In the latter case, the rationals must be
replaced by ratios of polynomials in $\K[X]$. We will prove that any
set $E$ invariant under rational translations in either of these
situations must also satisfy $\cH^f (E \cap F) = 0$ or $+ \infty$ for any
dimension function $f$ and any set $F$, which is the closure of an
open set of positive Haar measure.

\vskip 6mm

\centerline{\bf 2. Statement of the results}

\vskip 6mm

Our main result provides the appropriate analogue of Theorem A for
translation invariant outer measures (in particular Hausdorff
measures) as well as for higher dimensions. It rests on a clever idea
of Jarn\'\i k [10,11]. 

\proclaim Theorem 1. Let $E$ be a measurable subset of $\R^n$ which
does not have full Lebesgue measure. Assume that there exist integers
$\ell_1, \ldots, \ell_n \ge 2$ such that $E$ is invariant under any
translation by a vector of the form $(a_1/\ell_1^k, \dots,
a_n/\ell_n^k)$, with $\a = (a_1, \dots, a_n) \in \Z^n$, $k \in \Z$ and
$k \ge 1$. Let $\cM$ be a translation invariant outer measure and let
$F \subseteq \R^n$ be the closure of an open set of positive Lebesgue
measure. Then either $\cM(E \cap F) = 0$ or $+ \infty$.

Clearly, this theorem looks somewhat weaker than Theorem A, as it does
not show that such a set $E$ is either null or full with respect to
$\cM$. However, by taking $\cM$ to be the Lebesgue measure in the
proof and making the obvious modifications, a simple proof of Theorem
A, independent of the Lebesgue density theorem, is obtained.

As an immediate application of Theorem 1, we answer a question posed
by Mauldin ([6], p. 231):

{\it Does there exist a gauge function $g$ such that the Hausdorff $g$ 
measure of the set of the Liouville numbers is positive and finite?}
 
Note that Mauldin speaks of gauge functions and not dimension
functions. We will use the latter terminology. Recall that the set
$\cL$ of Liouville numbers is 
$$
\cL := \biggl\{\xi \in \R : \hbox{for all $w>0$, there exists} \ens
{p \over q} \ens \hbox{such that} \ens 
\biggl| \xi - {p \over q} \biggr| < {1 \over q^w} \biggr\}.
$$ It is readily verified that given any Liouville number $\xi$ and
any rational number $a/b$, the real number $\xi + a/b$ is a Liouville
number.  This implies 
\proclaim Corollary 1. Let $I \subseteq \R$ be a non-empty open
interval. There does not exist a dimension function $f$ such that $0 <
\cH^f (\cL \cap I) < + \infty$. 

By Corollary 1, $\cH^f (\cL)$ is either $0$ or $+ \infty$, but as far as
we are aware, given a dimension function $f$, it is an open problem to
determine which of these is the value of $\cH^f (\cL)$.  There are
some partial results on this question: using a covering argument, it
is easy to prove $\cH^f (\cL) = 0$ if there exists $\de > 0$ such that
$\lim_{x \to 0} \, x^{-\de} f(x) = 0$.  Conversely, Baker [2] showed
that if $f$ satisfies $\liminf_{x \to 0} \, x^{-\de} f(x) > 0$ for any
$\de > 0$, then we have $\cH^f (\cL) = + \infty$.

In dimensions strictly greater than one, sets of vectors analogous to
the sets $\cK(\tau)$, $\BA$ and $\cL$ also exist (see {\it e.g.}\ [21]
for examples). As these are also invariant under rational
translations, Theorem 1 immediately implies the same conclusion for
these sets. 

We briefly mention other applications of Theorem 1 to Diophantine
approximation.  In order to classify the real numbers according to
their properties of algebraic approximation, Mahler [16] and Koksma
[13] introduced, for any positive integer $n$, the functions $w_n$ and
$w^*_n$ defined as follows. Let $\xi$ be a real number. We denote by
$w_n (\xi)$ the supremum of the real numbers $w$ for which there exist
infinitely many integer polynomials $P(X)$ of degree at most $n$
satisfying
$$
0 < |P(\xi)| \le \HH(P)^{-w}, 
$$
where $\HH(P)$ is the na\"\i ve height of $P(X)$, that is, the maximum
of the absolute values of its coefficients. Further, we denote by
$w_n^*(\xi)$ the supremum of the real numbers $w^*$ for which there
exist infinitely many real algebraic numbers $\al$ of degree at most
$n$ satisfying
$$
0 < |\xi - \al| \le \HH(\al)^{-w^*-1}, 
$$
where $\HH(\al)$ is the na\"\i ve height of $\al$, that is, the height
of its minimal polynomial over the integers.  For results on the
functions $w_n$ and $w_n^*$, the reader is referred to [5,21]. We
observe that these functions are invariant by rational translations,
thus, by Theorem A, for any positive real numbers $w$ and $w^*$, any
of the sets
$$
\{\xi \in \R : w_n(\xi) = w\}, \qquad
\{\xi \in \R : w_n^*(\xi) = w^*\}
$$ have either Lebesgue measure zero, or have full Lebesgue measure.
This result (first observed by Sprind\v zuk [20]) can be extended
to Hausdorff measure (or indeed any translation invariant outer
measure) thanks to Theorem 1.

\proclaim Corollary 2. Let $w$ and $w^*$ be in $(0, +\infty]$.
Let $n$ be a positive integer.  There does not exist a dimension
function $f$ such that the intersection of any of the sets
$$
\{\xi \in \R : w_n(\xi) = w\}, \qquad
\{\xi \in \R : w_n^*(\xi) = w^*\}
$$
with any closed non-trivial interval has positive and finite Hausdorff
$\cH^f$-measure.  

Theorem 1 also applies to the second example considered by Olsen [17],
namely the Besicovitch--Eggleston set $B({\bf p})$ of non-normal
numbers in a given integer basis $N \ge 2$ (see [17] for the
definition), since this set is invariant under translation by any
rational number $a/N^k$, with $a, k \in \Z$ and $k\ge 1$.  It then
follows that we have either $\cH^f (B({\bf p}) \cap I) = 0$ or $\cH^f
(B({\bf p}) \cap I) = + \infty$, for any dimension function $f$ and
any closed non-trivial interval $I$. This improves Theorem 3 of [17]
and an earlier result of Smorodinsky [19]. For a restricted class of
Hausdorff measures, the same conclusion has been shown to hold by Ma
and Wen [15]. By Theorem 1, the same conclusion also holds for the
Cartesian product of such sets $B_{N_1}({\bf p}_1) \times \cdots
\times B_{N_n}({\bf p}_n)$, where the $N_i$ denote the bases and the
${\bf p}_i$ denote the required distribution of digits. This requires
the full force of the theorem.

For the $p$-adic fields and the fields of formal power series over
a finite field, we prove the following result.

\proclaim Theorem 2. Let $V$ be an $n$-dimensional vector space over
either a $p$-adic field or a field of formal Laurent series with
coefficients from a finite field $\K$. Let $E \subseteq V$ be a Haar
measurable set which does not have full Haar measure. Assume that $E$
is invariant under any translation by any vector of rationals
$(p_1/q_1, \ldots, p_n/q_n)$ in the $p$-adic case or any vector of
ratios $(p_1/q_1, \ldots, p_n/q_n)$ where $p_i,q_i \in \K[X]$, the
polynomial ring over $\K$, with the $q_i \neq 0$ in the case of formal
power series. Let $\cM$ be a translation invariant outer measure on
$V$ and let $F \subseteq V$ be the closure of an open set of positive
Haar measure. We then have $\cM(E \cap F) = 0$ or $\cM(E \cap F) = +
\infty$.

Note that Theorem 2 immediately implies that there is no dimension
function such that the sets analoguous to the set of Liouville numbers
have positive and finite measure. It has previously been shown that
the ordinary Hausdorff dimension of these sets is zero [1,14]. A full
analogue of Theorem 1 is possible in the case of formal power series,
which implies the same conclusion for the analogues of the
Besicovitch--Eggleston sets. However, the proof of the present
Theorem~2 is more elegant, and the reader should have no trouble
filling in the details to prove the full analogue of Theorem 1.
 
\vskip 6mm

\centerline{\bf 3. An important lemma}

\vskip 6mm

The fundamental tool in the proofs is the following lemma. In the case
of the real numbers, it is implicit in Jarn\'\i k's papers [10,11]. A
certain weak form of `quasi-independence' with respect to an outer
measure $\cal {M}$ can be defined; and it implies a `0-$\infty$' law
for $\cal {M}$ and so for any Hausdorff measure. We prove the lemma in
high generality.

\proclaim Lemma 1.  Let $F \subset \G$, where $\G$ is a locally
compact group and $F$ has finite Haar measure. Let $\mu$ denote the
restriction of the Haar measure on $\G$ to $F$, normalised so that
$\mu(F)=1$. Let $E \subset F$ be a measurable set with $\mu(E) < 1$
and let $\cal{M}$ be an outer measure on $F$. Suppose that for every
open ball $B(c,\rho) = \{x \in F : d(x,c)<\rho\} \subset F$,
$$ 
{\cal{M}}(E \cap B(c,\rho)) \, \le \, \mu(B(c,\rho)) {\cal{M}}(E).
$$
Then the measure ${\cal{M}}(E)$ of $E$ is either 0 or infinity.

Note that dividing both sides of the inequality in Lemma 1 yields
$$
{{\cal{M}}(E \cap B(c,\rho))\over {\cal{M}} (B(c,\rho))}
 \, \le \, {\mu(B(c,\rho))\over {\cal{M}}(B(c,\rho))} {\cal{M}}(E).
$$
The left hand side suggests a `quasi-density' (though of course the
limit need not exist) or a `conditional measure', while the right hand
side might be related to a notion analogous to `absolute continuity'
of the measures (although such a notion may not be appropriate in the
context of outer measures).

\vskip 2mm

\noindent {\bf Proof.} Assume the contrary, {\it i.e.}, assume $0 <
{\cal{M}}(E) < +\infty $.  Since $\mu(E)<1$, there exists a cover of $E$
by open balls $B(c_j, \rho_j)$ such that
$$
\sum_j \mu(B(c_j, \rho_j)) < 1\, .
$$
By assumption, 
$$
\eqalign{
0 < {\cal{M}}(E) &= {\cal{M}}\, \biggl(E \cap \biggl(\bigcup_j B(c_j,
\rho_j) \biggr) \biggr) = {\cal{M}}\, \biggl(\bigcup_j
(E \cap B(c_j, \rho_j) )\biggr) \cr
& \le \sum_j {\cal{M}} \biggl(E \cap B(c_j, \rho_j)\biggr)  \le \,
{\cal{M}}(E)\sum_j \mu(B(c_j, \rho_j)) < {\cal{M}}(E).
}
$$
This gives the desired contradiction.
\cqfd

\medskip

\vskip 6mm

\centerline{\bf 4. Proof of Theorem 1}

\vskip 6mm

It suffices to prove the theorem in the case when $F =
[0,h_1/\ell_1^{k'}] \times \cdots \times [0,h_n/\ell_n^{k'}]$. Indeed,
suppose that $F'$ is the closure of an open set. Then there is a
vector $\r$ of the form from Theorem 1 and an $F$ such that $F+\r
\subseteq F'$. If $\cM(E \cap F) = +\infty$, the translation
invariance of $E$ implies that $\cM(E \cap F') = +\infty$. 

Suppose on the other hand that $\cM(E \cap F) = 0$. Then we may cover
$F'$ by countably many translates of the form from Theorem 1. In this
case, translation invariance of $E$ implies that $\cM(E \cap F) = 0$.
We will prove the theorem in the case when $F = [0,1]^n$, the closed
unit hypercube. This is to avoid the notational complications of
additional indices. The reader should have no trouble filling in the
details for other hyperboxes of the above form.

Let $\j = (j_1, \dots, j_n) \in \Z^n$ and let
$$
I (\j;k) = \biggl[ {j_1 \over \ell_1^k}, {j_1+1 \over \ell_1^k}
\biggr] \times \cdots \times \biggl[ {j_n \over \ell_n^k}, {j_n+1 \over
\ell_n^k} \biggr]. 
$$
As $E$ is assumed to be invariant under translations by vectors of the
form $(j_1/\ell_1^k, \ldots, j_n/\ell_n^k)$, we have
$$
E \cap I(\j;k) = \biggl( E \cap \biggl[0, {1 \over \ell_1^k} \biggr]
\times \cdots \times \biggl[0, {1 \over \ell_n^k} \biggr] \biggr)
+ (j_1/\ell_1^k, \ldots, j_n/\ell_n^k).
$$ 
Since the outer measure $\cM$ is assumed to be translation invariant,
we get 
$$
\eqalign{
\cM ( E \cap [0, 1]^n) &= \sum_{j_1 = 0}^{\ell_1^k - 1} \, \cdots \,
\sum_{j_n = 0}^{\ell_n^k - 1} \, \cM \bigl( E \cap I(\j;k)
\bigr) \cr 
&= \ell_1^k \cdots \ell_n^k \, \cM \biggl( E \cap \biggl[0, {1 \over
\ell_1^k} \biggr] \times \cdots \times \biggl[0, {1 \over \ell_n^k}
\biggr] \biggr).
}
$$
Consequently, for any $\j \in \Z^n$ and any $k$ with $k \ge 1$ and $0
\le j_i \le \ell_i^k - 1$, we have
$$
\cM \bigl( E \cap I(\j;k) \bigr) = {1 \over \ell_1^{k} \cdots
\ell_n^k} \, \cM ( E \cap [0, 1]^n).  \eqno (1)
$$

We endow $\R^n$ with the metric induced by the norm $|x|_\infty =
\max\{|x_1|, \ldots, |x_n|\}$. In this metric, an open ball
$B(c,\rho)$ is the Cartesian product of $n$ open intervals, {\it
i.e.}, $B(c,\rho) = (a_1, b_1) \times \cdots \times (a_n, b_n)$. 
Considering the expansions of its endpoints in some base $\ell$, every 
real open interval $(a, b)$ can be represented as a union of a
countable set of intervals $[j/\ell^k, (j+1)/\ell^k]$ in such a way
that
$$
(a, b) = \bigcup_{j, k} \, \biggl[ {j \over \ell^k}, {j+1 \over \ell^k} 
\biggr] \qquad {\rm and} \qquad b-a = \sum_{j, k} \, {1 \over \ell^k}.
$$
We do this for each coordinate, where we expand the $i$'th coordinate
interval in base $\ell_i$. If necessary, we subdivide intervals again to
obtain a representation for $B(c,\rho)$ such that  
$$
B(c,\rho) = \bigcup_{\j, k} \, I(\j;k) \qquad {\rm and} \qquad
\mu\bigl(B(c,\rho)\bigr) = \sum_{\j, k} \, {1 \over \ell_1^k \cdots
\ell_n^k},
$$
where $\mu$ denotes the Lebesgue measure on $\R^n$. Hence, we get
$$
\cM \bigl(B(c,\rho)\bigr) = \cM \, \biggl( \bigcup_{\j, k} \,
I(\j;k) \biggr),
$$
and, since $\cM ( \cdot)$ is an outer measure, it follows from (1)
that
$$ 
\eqalign{ 
\cM \bigl(B(c,\rho) \cap E \cap [0, 1]^n\bigr) & = \cM
\biggl( \bigcup_{\j, k} \, E \cap I(\j;k) \cap [0, 1]^n \biggr) \cr
& \le \sum_{\j, k} \, \cM \bigl( E \cap I(\j;k) \cap [0, 1]^n \bigr)
\cr 
& \le \cM \bigl(E \cap [0, 1]^n\bigr) \, \sum_{\j, k} \, {1 \over
\ell_1^k \cdots \ell_n^k} \cr
&= \mu\bigl(B(c,\rho)\bigr) \, \cM (E \cap
[0, 1]^n).} \eqno (2)
$$

Now by hypothesis, E does not have full measure and by taking
$X=[0,1]^n$ in Lemma~1 (as we can without loss of generality), we can
suppose that $\mu(E\cap [0,1]^n)<1$.  
It follows from Lemma 1 that $\cM(E \cap [0,1]^n) = 0$ or $\cM(E \cap
[0,1]^n) = + \infty$. This is the statement of Theorem 1 in the
case when $F = [0,1]^n$. For other sets $F$ the theorem follows
analogously.

Furthermore, taking $\cM$ to be the outer Lebesgue measure in the
proof, the set $E \cap [0,1]^n$ must have outer Lebesgue measure
either $0$ or $\infty$. Clearly, $\cM(E \cap [0,1]^n) \neq \infty$,
and so we have given a simple proof of Theorem A. \cqfd

\vskip 6mm

\centerline{\bf 5. Proof of Theorem 2}

\vskip 6mm

The proof relies on the same idea as the proof of Theorem 1 and is
almost identical for $p$-adics and formal power series.  As before, it
suffices to consider the case when $F = B(0,1)^n$, the unit hypercube
in $V$. We let $\r$ denote some vector with rational coordinates in
the $p$-adic case and ratios of polynomials in the case of formal
power series. Let $k$ denote $p$ in the $p$-adic case and $\vert \K
\vert$ in the case of formal power series.

Let $B(c,\rho)$ denote a closed ball centred at $c$ with radius $\rho$
in the metric induced by the height $\max\{|x_1|, \ldots, |x_n|\}$,
where $|\cdot|$ denotes the absolute value on the base field. Note
that because of the definition of the metric of the underlying space,
for any $\rho > 0$, we have $B(c,\rho) = B(c, k^{-r})$ for some $r \in \Z$. We
may therefore restrict ourselves to considering balls with radii of
this form.

As in the real case, translation invariance implies that for any $\r$
and any $r \in \Z$,
$$
E \cap B(\r, k^{-r}) = E \cap B(0, k^{-r}) + \r.
$$ Furthermore, using the ultrametric property of the underlying
spaces, it is possible to tile the unit ball $B(0,1)$ with $k^{nr}$
disjoint balls of radius $k^{-r}$. Using the translation invariance of
$\cM$, we get
$$
\cM(E \cap B(0,1)) = k^{nr} \cM (E \cap B(0,k^{-r})).
$$
Hence, for any ball of radius $k^{-r}$ centred at $\r$, we have
$$
\cM ( E \cap B(\r, k^{-r}))
= k^{-nr} \, \cM ( E \cap B(0,1)).  \eqno (3)
$$

Now, the spaces considered are ultrametric, and so any interior point
of a ball may be taken as the centre of the ball. Furthermore, the set
of elements $\r$ is dense in the spaces by construction, so any ball
of positive radius has such a point as an interior point. Hence, for
any ball $B(c, k^{-r})$, there is an $\r$ such that $B(c, k^{-r}) =
B(\r, k^{-r})$. Using (3), for any ball $B(c, k^{-r})$, we obtain
$$
\eqalign{
\cM(E \cap B(c, k^{-r}) \cap B(0,1)) &= \cM(E \cap B(\r, k^{-r})
\cap B(0,1)) \cr
& \leq k^{-nr} \cM(E \cap B(0,1)) = \mu(B(c,k^{-r})) \cM(E \cap
B(0,1)),
}
$$ 
where $\mu$ is the Haar measure, normalised so that the closed unit
ball has measure $1$. On supposing that $\mu(E \cap B(0,1)) < 1$, we
may invoke Lemma 1 to prove the theorem. In the case when $\mu(E \cap
B(0,1)) = 1$, considering a union of translates of this set gives an
analogue of Theorem A.  \cqfd

\vskip 10mm

\centerline{\bf References}

\vskip 8mm

\item{[1]} A. Abercrombie,
{\it The Hausdorff dimension of some exceptional sets of $p$-adic
integer matrices}, J. Number Theory 53 (1995), 311--341. 
\sm
\item{[2]} R. C. Baker,
{\it On approximation with algebraic numbers of bounded degree},
Mathematika 23 (1976), 18--31.
\sm
\item{[3]} V. Beresnevich, D. Dickinson and S. Velani,
{\it Measure theoretic laws for lim sup sets}, preprint.
\sm
\item{[4]} Y. Bugeaud, 
{\it Approximation par des nombres alg\'ebriques de degr\'e born\'e
et dimension de Hausdorff}, 
J. Number Theory 96 (2002), 174--200.
\sm
\item{[5]} Y. Bugeaud, Approximation by algebraic numbers,
Cambridge Tracts in Mathematics, Cambridge University Press,
Cambridge, 2004.
\sm
\item{[6]} M. Cs\"ornyei, {\it Open Problems. Compiled and edited by
Marianne Cs\"ornyei}, Proceedings from the conference `Dimensions and
Dynamics', Miskolc, Hungary, July 20--24, 1998. Periodica Math.
Hung. 37 (1998), 227--237.
\sm
\item{[7]} M. M. Dodson and S. Kristensen,
{\it Hausdorff dimension and Diophantine approximation},
`Fractal Geometry and Applications: A Jubilee of Benoit Mandelbrot', 
Proceedings of Symposia in Pure Mathematics, American Mathematical
Society.
To appear.
\sm
\item{[8]} K. Falconer, 
The geometry of fractal sets, 
Cambridge Tracts in Mathematics 85, Cambridge University Press, 1985.
\sm
\item{[9]}  V. Jarn\'\i k, 
{\it Diophantischen Approximationen und Hausdorffsches Mass},
Mat. Sbornik 36 (1929), 371--382.
\sm
\item{[10]}  V. Jarn\'\i k, 
{\it N\v ekolik pozn\'amek o Hausdorffov\v e m\'\i \v re},
Rozpravy T\'r. \v Cesk\'e Akad 40, c. 9 (1930), 8 p.
\sm
\item{[11]} V. Jarn\'\i k, 
{\it Quelques remarques sur la mesure de M. Hausdorff},
Bull. Int. Acad. Sci. Boh\^eme (1930), 1--6.
\sm
\item{[12]}  V. Jarn\'\i k, 
{\it \"{U}ber die simultanen diophantischen Approximationen},
Math. Z. 33 (1931), 505--543.
\sm
\item{[13]} J. F. Koksma,
{\it \"Uber die Mahlersche Klasseneinteilung der transzendenten Zahlen
und die Approximation komplexer Zahlen durch algebraische Zahlen},
Monats. Math. Phys. 48 (1939), 176--189.
\sm  
\item{[14]} S. Kristensen,
{\it On well-approximable matrices over a field of formal series},
Math. Proc. Cambridge Philos. Soc. 135 (2003), 255--268.
\sm
\item{[15]} J. Ma and Z. Wen,
{\it Hausdorff and packing measure of sets of generic points: A
zero-infinity law}, J. London Math. Soc. (2) 69 (2004), 383--406.
\sm
\item{[16]} K. Mahler,
{\it Zur Approximation der Exponentialfunktionen und des
Logarithmus. I, II},
J. reine angew. Math. 166 (1932), 118--150.
\sm
\item{[17]} L. Olsen, {\it On the dimensionlessness of invariant sets},
Glasgow Math. J. 45 (2003), 539--543.
\sm
\item{[18]} C. A. Rogers, Hausdorff Measures, Cambridge University
Press, Cambridge, 1970.
\sm
\item{[19]} M. Smorodinsky, {\it Singular measures and Hausdorff
measures}, Israel J. Math. 7 (1969), 203--206.
\sm
\item{[20]} V. G. Sprind\v zuk,
{\it On some general problems of approximating numbers by algebraic numbers},
Litovsk. Mat. Sb. 2 (1962), 129--145 (in Russian).

\item{[21]} V. G. Sprind\v zuk,
Metric theory of Diophantine approximation,
Izdat. ``Nauka'', Moscow, 1977 (in Russian). English translation by
R. A. Silverman, Scripta Series in Mathematics,
John Wiley \& Sons, New-York--Toronto--London, 1979.

\vskip 5mm

\noi Y. Bugeaud

\noi Universit\'e Louis Pasteur

\noi U. F. R. de math\'ematiques

\noi 7, rue Ren\'e Descartes

\noi 67084 STRASBOURG Cedex

\vskip 2mm

\noi e-mail : {\tt bugeaud@math.u-strasbg.fr}

\vskip 5mm

\noi M. M. Dodson

\noi Department of Mathematics

\noi University of York

\noi Heslington

\noi York

\noi YO10 5DD

\noi UK

\vskip 2mm

\noi e-mail : {\tt mmd1@york.ac.uk}

\vskip 5mm

\noi S. Kristensen

\noi School of Mathematics

\noi University of Edinburgh 

\noi  James Clerk Maxwell Building

\noi  King's Buildings

\noi  Mayfield Road 

\noi  Edinburgh 

\noi  EH9 3JZ 

\noi  UK

\vskip 2mm

\noi e-mail : {\tt Simon.Kristensen@ed.ac.uk}

\bye